\documentclass[letterpaper,conference]{ieeeconf}
\IEEEoverridecommandlockouts                              
\usepackage[utf8]{inputenc}
\usepackage{textcomp}

\usepackage{graphicx}
\usepackage{amsmath}
\usepackage[version=4]{mhchem}
\usepackage{siunitx}
\usepackage{graphicx}
\usepackage{longtable,tabularx}
\usepackage[style=ieee ,sorting=none,maxbibnames=99,giveninits, doi=false]{biblatex}
\AtBeginBibliography{\small}

\let\labelindent\relax
\usepackage{amsmath} 
\usepackage{amsfonts}
\usepackage{amsthm}

\theoremstyle{definition}

\newtheorem{remark}{Remark}

\usepackage{float} 
\usepackage[skip=0em,font=footnotesize]{caption}

\usepackage[margin = 1in]{geometry}

\usepackage{tikz}
\usetikzlibrary{shapes,arrows,calc,positioning}
\tikzstyle{bigblock} = [draw, fill=blue!20, rectangle, 
    minimum height=6em, minimum width=8em]
\tikzstyle{medblock} = [draw, fill=blue!20, rectangle, 
    minimum height=4em, minimum width=4em]    
\tikzstyle{mux} = [draw, fill=black!20, rectangle, 
    minimum height=5em, minimum width=0.1em]    
\tikzstyle{smallblock} = [draw, fill=blue!20, rectangle, 
    minimum height=2em, minimum width=3em]
    
\tikzstyle{data_block} = [draw, fill=green!20, rectangle, 
    minimum height=2em, minimum width=3em]
\tikzstyle{ops_block} = [draw, fill=blue!20, rectangle, 
    minimum height=2em, minimum width=3em]    
\tikzstyle{est_block} = [draw, fill=red!20, rectangle, 
    minimum height=2em, minimum width=3em]    
    
\tikzstyle{sum} = [draw, fill=blue!20, circle, node distance=1cm,minimum height=0.5cm]
\tikzstyle{signal} = [coordinate]
\tikzstyle{pinstyle} = [pin edge={to-,thin,black}]
\tikzstyle{block} = [draw, fill=blue!20, rectangle, 
    minimum height=3em, minimum width=9em]
\tikzstyle{blockS} = [draw, fill=blue!20, rectangle, 
    minimum height=3em, minimum width=4em]    
\tikzstyle{input} = [coordinate]
\tikzstyle{output} = [coordinate]
\usetikzlibrary{matrix}

\usetikzlibrary{positioning}
\usetikzlibrary{math}

\usepackage{hyperref}
\usepackage{xcolor}
\hypersetup{
    colorlinks,
    linkcolor={blue!100!black},
    citecolor={blue!50!black},
    urlcolor={blue!80!black}
}

\newcommand{\bc}{\begin{center}}
\newcommand{\ec}{\end{center}}
\newcommand{\benum}{\begin{enumerate}}
\newcommand{\eenum}{\end{enumerate}}

\newcommand{\matl}{\left[ \begin{array}}
\newcommand{\matr}{\end{array} \right]}

\renewcommand{\matl}{\begin{bmatrix}}
\renewcommand{\matr}{\end{bmatrix}}

\newcommand{\matls}{\left[ \begin{smallmatrix}}
\newcommand{\matrs}{\end{smallmatrix} \right]}
\newcommand{\isdef}{\stackrel{\triangle}{=}}

\newcommand{\rmC}{{\rm C}}

\newcommand{\rmS}{{\rm S}}

\newcommand{\rmc}{{\rm c}}

\newcommand{\rms}{{\rm s}}

\newcommand{\SG}{{\mathcal G}}

\usepackage{enumitem,amssymb}
\newlist{todolist}{itemize}{2}
\setlist[todolist]{label=$\square$}
\usepackage{pifont}

\usepackage{booktabs} 
\title{Reference Governor for Constraint Enforcement in Data-Driven Model Adaptive Control}
\title{Constraint Enforcement via Reference Governor in \\ Indirect Adaptive Control}
\title{Constraint Enforcement via Reference Governor in \\ Dynamic Mode Adaptive Control}

\author{
Parham Oveissi and
Ankit Goel
\thanks{Parham Oveissi is a graduate student in the Department of Mechanical Engineering, University of Maryland, Baltimore County, 1000 Hilltop Circle, Baltimore, MD 21250. {\tt\small parhamo1@umbc.edu}}%
\thanks{Ankit Goel is an Assistant Professor in the Department of Mechanical Engineering, University of Maryland, Baltimore County,1000 Hilltop Circle, Baltimore, MD 21250. {\tt\small ankgoel@umbc.edu }}%
}

\begin{document}

\maketitle

\begin{abstract}
This paper integrates reference governance with Dynamic Mode Adaptive Control (DMAC) to enforce constraints without a prior plant model.
DMAC identifies and updates a discrete-time model online from measured state and input data, and the reference governor (RG) uses this identified model, rather than the true plant dynamics, to predict constrained closed-loop behavior.
Because the prediction model changes as DMAC adapts, a governed reference that was previously feasible may become infeasible, and the classical interpolation-based RG may be unable to move the reference in the direction required to recover predicted constraint satisfaction.
The proposed DMAC--RG framework instead searches directly over the governed reference rather than over an interpolation factor.
When feasible references exist, the governed reference is selected to track the desired command, with an optional penalty on deviation from the previously applied reference.
When the finite-horizon feasible set is empty, the reference is selected to minimize predicted constraint violation while retaining the same tracking preference.
This work is a numerical investigation and does not provide theoretical guarantees of recursive feasibility, constraint satisfaction, or closed-loop stability.
Numerical examples on a linear system, a nonlinear Van der Pol oscillator, and a high-dimensional Burgers equation under partial-state measurement demonstrate that the proposed RG enforces constraints using only the model identified online by DMAC, without access to the true plant dynamics or the complete system state.
\end{abstract}

\textit{\bf keywords:} Constrained control, Adaptive control, Indirect adaptive control.

\section{INTRODUCTION}

Constraints on states, outputs, and control inputs arise naturally in
control systems because of actuator limitations, operating envelopes,
safety requirements, and physical restrictions.
A controller designed primarily for stabilization or tracking may produce
satisfactory nominal closed-loop behavior while still violating such
constraints in response to demanding commands.
Consequently, constraint management is an important consideration when
applying a controller to a physical system.

Reference governors (RGs) provide an attractive framework for incorporating
constraint handling into an existing closed-loop system.
Rather than redesigning the underlying controller, an RG acts as an
add-on supervisory element that modifies the commanded reference whenever
application of the original command is predicted to violate prescribed
constraints
\cite{gilbert1994nonlinear,gilbert1995discrete,garone2017reference}.
When the desired command is admissible, it is passed to the controller;
otherwise, the governed reference is modified to maintain admissibility
while remaining close to the desired command.
This separation between the nominal controller and the constraint-management
mechanism has led to RG formulations for linear, nonlinear, uncertain,
distributed, and networked systems
\cite{gilbert2002nonlinear,garone2017reference}.

A central requirement of reference governance is a model with which the
future response of the constrained variables can be predicted.
Classical RG formulations are therefore naturally suited to systems for
which an adequate plant or closed-loop model is available.
When the dynamics are uncertain, robust RG formulations can account for
bounded uncertainty, although the resulting constraint-admissible sets may
be conservative.
Adaptive extensions have consequently been developed to exploit online
information about uncertain parameters.
For example, parameter-adaptive RG methods use online parameter estimates
together with robust constraint-admissible sets
\cite{chakrabarty2023parameter}, while integrated adaptive-control and RG
architectures have used uncertainty compensation and constraint tightening
to provide constraint handling in the presence of state-dependent
uncertainties \cite{zhao2023integrated}.
Data-driven approaches have also been investigated, including learning
invariant sets directly from data and constructing predictive models for
unknown nonlinear dynamics \cite{kashani2023data}.

The setting considered in this work is motivated by applications in which
an accurate plant model may not be available, but constraint management is
nevertheless required.
One motivating application is the control of shock-train dynamics in
supersonic propulsion systems.
Previous work has demonstrated the use of adaptive and data-driven control
for thrust regulation in solid-fuel ramjets, including configurations with
variable-geometry inlets and in-situ thrust estimation.
These applications motivate the development of constraint-management
mechanisms that can operate using dynamics identified directly from
measured data rather than relying on a predetermined plant model \cite{oveissi2026adaptive, oveissi2026thrust}.

Dynamic Mode Adaptive Control (DMAC) provides such an adaptive-control
framework.
DMAC identifies a discrete-time linear model from measured input--state
data and uses the identified dynamics for online controller synthesis
\cite{oveissi2026model}.
The model available to the controller therefore evolves as new data are
acquired.
The objective of this work is to use the same online-identified dynamics
for RG prediction, thereby incorporating reference governance into DMAC
without requiring a separately specified plant model.

This integration introduces a difficulty because the prediction model used
by the RG changes as DMAC updates its identified dynamics.
Consequently, a governed reference that was feasible under the model
available at one sample need not remain feasible under the updated model.
The classical interpolation-based RG restricts the new governed reference
to the segment between the previously applied reference and the desired
command.
If the previous reference loses feasibility, this restriction may prevent
the governor from moving the reference in the direction required to
recover predicted constraint satisfaction or reduce predicted constraint
violation.
In addition, the finite-horizon feasible reference set may become empty,
requiring a reference-selection rule even when no candidate satisfies all
predicted constraints.

The main contribution of this paper is the integration of reference
governance with DMAC using the dynamics identified online by the adaptive
controller.
Rather than parameterizing the governed reference through the classical
interpolation factor, the proposed RG searches directly over the governed
reference.
When feasible candidates exist, the governed reference is selected by
balancing tracking of the desired command with variation from the
previously applied reference.
When the finite-horizon feasible set is empty, the reference is selected
to reduce the accumulated predicted output and input constraint violation
while retaining the same reference-tracking preference.
The resulting formulation allows the governed reference to move outside
the classical interpolation segment and provides a defined
reference-selection rule when the predicted feasible set is empty.

The present work is intended as a preliminary numerical investigation of
the integration of reference governance with DMAC.
The proposed framework is evaluated entirely through numerical experiments,
and no theoretical guarantees of recursive feasibility, constraint
satisfaction, or closed-loop stability are established.
The numerical studies are used to investigate whether the online-identified
DMAC dynamics can provide a useful prediction model for reference
governance across systems of differing complexity.
A theoretical treatment of the resulting time-varying prediction model
and associated feasibility and stability properties is left for future
work.

The remainder of the paper is organized as follows.
Section~\ref{sec:DMAC} briefly reviews the DMAC framework and the
online-identified dynamics used for RG prediction.
Section~\ref{sec:classicalRG} reviews the classical reference governor
and discusses its limitations when the prediction model changes online.
Section~\ref{sec:Unrestricted_RG} presents the proposed unrestricted
reference-search formulation, including reference selection when the
finite-horizon feasible set is empty.
Section~\ref{sec:numerical_results} presents numerical results for linear
and nonlinear systems, including a high-dimensional Burgers equation
example with partial-state feedback.
Finally, Section~\ref{sec:conclusion} summarizes the results.

\section{DYNAMIC MODE ADAPTIVE CONTROL}
\label{sec:DMAC}

Dynamic Mode Adaptive Control (DMAC) provides the baseline adaptive-control architecture considered in this work \cite{oveissi2026model}.
Let $x_k$ denote the system state and let $\xi_k\in\mathbb{R}^{n_\xi}$ denote the measured portion of the state available to DMAC.
At each sample $k$, DMAC identifies the local discrete-time dynamics
\begin{align}
    \xi_{k+1}
    &=
    \hat A_k \xi_k
    +
    \hat B_k u_k,
    \label{eq:dmac_identified_model}
\end{align}
from measured state and input data and updates the corresponding feedback controller online.

In this work, the regulated output is
\begin{align}
    y_k
    =
    C\xi_k,
    \label{eq:dmac_output}
\end{align}
and integral action is introduced according to
\begin{align}
    q_{k+1}
    =
    q_k
    +
    T_s
    \left(
        g_k-C\xi_k
    \right),
    \label{eq:dmac_integrator}
\end{align}
where $g_k$ is the governed reference.
The feedback law is written as
\begin{align}
    u_k
    =
    K_{\xi,k}\xi_k
    +
    K_{q,k}q_k.
    \label{eq:dmac_feedback}
\end{align}
The identified matrices $\hat A_k$ and $\hat B_k$ and the corresponding feedback gains are updated online by DMAC.
Further details of the identification and adaptive-control formulation are given in \cite{oveissi2026model}.

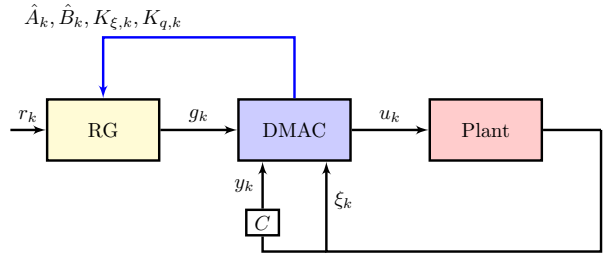
\begin{figure}[ht]
    \centering
    \resizebox{0.98\columnwidth}{!}{%
    \begin{tikzpicture}[
        >=latex',
        line width=1.2,
        block/.style={
            draw=black,
            rectangle,
            minimum width=5.5em,
            minimum height=3.0em,
            align=center
        }
    ]


    \node[
        block,
        fill=yellow!20
    ] (RG) at (0,0) {\textbf{RG}};

    \node[
        block,
        fill=blue!20
    ] (DMAC) at (3.3,0) {\textbf{DMAC}};

    \node[
        block,
        fill=red!20
    ] (Plant) at (6.6,0) {\textbf{Plant}};


    \draw[->]
        (-1.6,0)
        -- node[above] {$r_k$}
        (RG.west);

    \draw[->]
        (RG.east)
        -- node[above] {$g_k$}
        (DMAC.west);

    \draw[->]
        (DMAC.east)
        -- node[above] {$u_k$}
        (Plant.west);

    \draw[-]
        (Plant.east)
        -- ++(1.0,0)
        coordinate (out);


    \coordinate (xibranch) at
        ([xshift=0.55cm,yshift=-1.55cm]DMAC.south);

    \draw[->]
        (out)
        |- (xibranch)
        -- node[pos=0.58,right] {$\xi_k$}
        ([xshift=0.55cm]DMAC.south);


    \node[
        draw=black,
        rectangle,
        fill=white,
        minimum width=1.6em,
        minimum height=1.25em,
        inner sep=1pt
    ] (Cblock) at
        ([xshift=-0.55cm,yshift=-1.05cm]DMAC.south)
        {$C$};

    \draw[-]
        (xibranch)
        -| (Cblock.south);

    \draw[->]
        (Cblock.north)
        -- node[midway,left] {$y_k$}
        ([xshift=-0.55cm]DMAC.south);


    \draw[->,blue]
        (DMAC.north)
        -- ++(0,1.05)
        -|
        node[pos=0.5,above,black]
        {$\hat A_k,\hat B_k,K_{\xi,k},K_{q,k}$}
        (RG.north);

    \end{tikzpicture}%
    }
    \vspace{1em}
    \caption{DMAC--RG architecture using the online identified dynamics.}
    \label{fig:DMAC_RG_architecture}
\end{figure}

The contribution considered in this work is the integration of a reference governor with DMAC.
As shown in Fig.~\ref{fig:DMAC_RG_architecture}, the RG uses the dynamics identified online by DMAC to predict constraint satisfaction.
At each sample, the current estimates $\hat A_k$ and $\hat B_k$, together with the current feedback gains, are supplied to the RG.
The RG determines the governed reference $g_k$, which is then passed to DMAC as the reference command for the closed-loop system.
Thus, the same model identified online for adaptive control is also used by the RG for constraint prediction.

\section{CLASSICAL REFERENCE GOVERNOR}
\label{sec:classicalRG}

A reference governor (RG) is a supervisory constraint-handling mechanism
that modifies the reference supplied to an otherwise closed-loop system.
Rather than applying the desired reference $r_k$ directly, the RG
determines a governed reference $g_k$ such that the predicted closed-loop
response satisfies prescribed state and input constraints.
The controller itself is not modified by the RG.
Instead, constraint satisfaction is achieved by restricting the reference
command presented to the closed-loop system.

To determine whether a candidate governed reference is admissible, the RG
predicts the future closed-loop response assuming that the candidate
reference is held constant over the prediction horizon.
Consider the discrete-time closed-loop prediction model
\begin{align}
    z_{i+1}
    &=
    A z_i
    +
    B g,
    \label{eq:classical_system}
\end{align}
where $z_i\in\mathbb{R}^{n_z}$ is the closed-loop state and $g$ is the
constant candidate governed reference.
The state $z_i$ may include both plant and controller states required to
predict the closed-loop response.

The constrained plant and input quantities are obtained from the
closed-loop state according to
\begin{align}
    y_{\rmc,i} &= C_y z_i,
    \\
    u_{\rmc,i} &= C_u z_i,
\end{align}
and are required to satisfy
$y_{\rmc,\min}
\leq
y_{\rmc,i}
\leq
y_{\rmc,\max}$
and
$u_{\rmc,\min}
\leq
u_{\rmc,i}
\leq
u_{\rmc,\max}$.
A candidate reference $g$ is admissible if the trajectory predicted from
the current closed-loop state satisfies these constraints over the
prediction horizon.

The classical scalar RG searches for an admissible reference along the
line segment connecting the previously applied governed reference
$g_{k-1}$ and the desired reference $r_k$.
A candidate reference is parameterized as
\begin{align}
    g(\kappa)
    =
    g_{k-1}
    +
    \kappa
    \left(
        r_k-g_{k-1}
    \right),
    \label{eq:classical_interpolation}
\end{align}
where $\kappa\in[0,1].$
Thus, $\kappa=0$ retains the previously applied reference, whereas
$\kappa=1$ applies the desired reference.

For a finite prediction horizon $N_h$, the RG selects the largest
admissible interpolation factor according to
\begin{align}
    \kappa_k^\star
    &=
    \max_{\kappa\in[0,1]}
    \quad \kappa
    \label{eq:classical_optimization}
    \\
    \text{subject to} \qquad
    z_{i+1|k}
    &=
    A z_{i|k}
    +
    B g(\kappa),
    \nonumber\\
    y_{\rmc,\min}
    &\leq
    y_{\rmc,i|k}
    \leq
    y_{\rmc,\max},
    \nonumber\\
    u_{\rmc,\min}
    &\leq
    u_{\rmc,i|k}
    \leq
    u_{\rmc,\max},
    \nonumber
\end{align}
for $i=1,\ldots,N_h$.
The governed reference applied at time $k$ is then
\begin{align}
    g_k
    =
    g_{k-1}
    +
    \kappa_k^\star
    \left(
        r_k-g_{k-1}
    \right).
    \label{eq:kRG_result}
\end{align}
Hence, the classical RG moves the governed reference as far as possible
toward the desired reference while maintaining predicted constraint
satisfaction.

A central property of this formulation is that the search is restricted
to the line segment connecting $g_{k-1}$ and $r_k$.
Classical recursive-feasibility arguments rely on the previously applied
reference $g_{k-1}$ remaining admissible and therefore serving as a
fallback solution.
For a fixed prediction model and an appropriate invariant admissible set,
this property can be established because retaining $g_{k-1}$ continues a
previously admissible closed-loop trajectory
\cite{gilbert1995discrete,garone2017reference}.

This property is not guaranteed when the prediction model changes online.
In the DMAC setting, the model used by the RG is updated as new system
data become available.
Consequently, a reference that was admissible under the prediction model
at time $k-1$ may no longer be admissible under the model available at
time $k$.
In this case, choosing $\kappa=0$ does not necessarily provide a feasible
fallback reference.

The line-segment parameterization imposes an additional restriction.
Because $\kappa\in[0,1]$, the governed reference can move only between
$g_{k-1}$ and $r_k$.
If constraint satisfaction under the updated prediction model requires
the governed reference to move outside this interval, the classical RG
cannot generate such a command.
This limitation motivates the unrestricted reference-governor formulation
developed in Section~\ref{sec:Unrestricted_RG}.

\section{REFERENCE GOVERNOR WITH UNRESTRICTED REFERENCE SEARCH}
\label{sec:Unrestricted_RG}

To address the limitations identified in Section~\ref{sec:classicalRG},
the proposed RG searches directly over the governed reference $g$ rather
than restricting the search to the interpolation segment between
$g_{k-1}$ and $r_k$.
If feasible candidates exist, the reference closest to the desired
command is selected.
If the finite-horizon feasible set is empty, the reference is selected to
minimize the predicted constraint violation, with a secondary preference
for the desired reference.

\subsection{Closed-Loop Prediction Model}

At each sample $k$, DMAC provides the current identified dynamics
$\hat A_k$ and $\hat B_k$ together with the current feedback gains
$K_{\xi,k}$ and $K_{q,k}$.
For the RG prediction, these quantities are frozen at their current values
over the prediction horizon.
Accordingly, let
$A=\hat A_k$, $B=\hat B_k$, $K_\xi=K_{\xi,k}$, and
$K_q=K_{q,k}$.

Using the measured state $\xi_{i|k}$ available to DMAC and the integral
state $q_{i|k}$, define the augmented RG prediction state
\begin{align}
    z_{i|k}
    \isdef
    \begin{bmatrix}
        \xi_{i|k}\\
        q_{i|k}
    \end{bmatrix}.
\end{align}
Importantly, $z_{i|k}$ is the state of the RG prediction model and need
not contain the complete physical state of the plant.
In the DMAC architecture considered here, the prediction model is
constructed from the measured state used for online identification
together with the controller integral state.

For the integral-feedback structure considered here, the frozen
closed-loop prediction model is
\begin{align}
    z_{i+1|k}
    =
    \overline A z_{i|k}
    +
    \overline B g,
    \label{eq:rg_augmented_prediction}
\end{align}
for $i=0,\ldots,N_h-1$, where
\begin{align}
    \overline A
    &=
    \begin{bmatrix}
        A+BK_\xi & BK_q\\
        -T_sC    & 1
    \end{bmatrix},
    \quad 
    \overline B
    =
    \begin{bmatrix}
        0\\
        T_s
    \end{bmatrix}.
\end{align}
The prediction is initialized using
\begin{align}
    z_{0|k}
    =
    \begin{bmatrix}
        \xi_k\\
        q_k
    \end{bmatrix}.
\end{align}

The candidate reference $g$ is held constant over the prediction horizon.
Thus,
\begin{align}
    z_{i|k}(g)
    =
    \overline A^i z_{0|k}
    +
    \sum_{j=0}^{i-1}
    \overline A^j
    \overline B g.
    \label{eq:rg_closed_form_prediction}
\end{align}
Hence, for the frozen model, the predicted closed-loop state is affine in
the candidate reference $g$.

\begin{figure}[ht]
    \centering
    \includegraphics[width=\columnwidth]
    {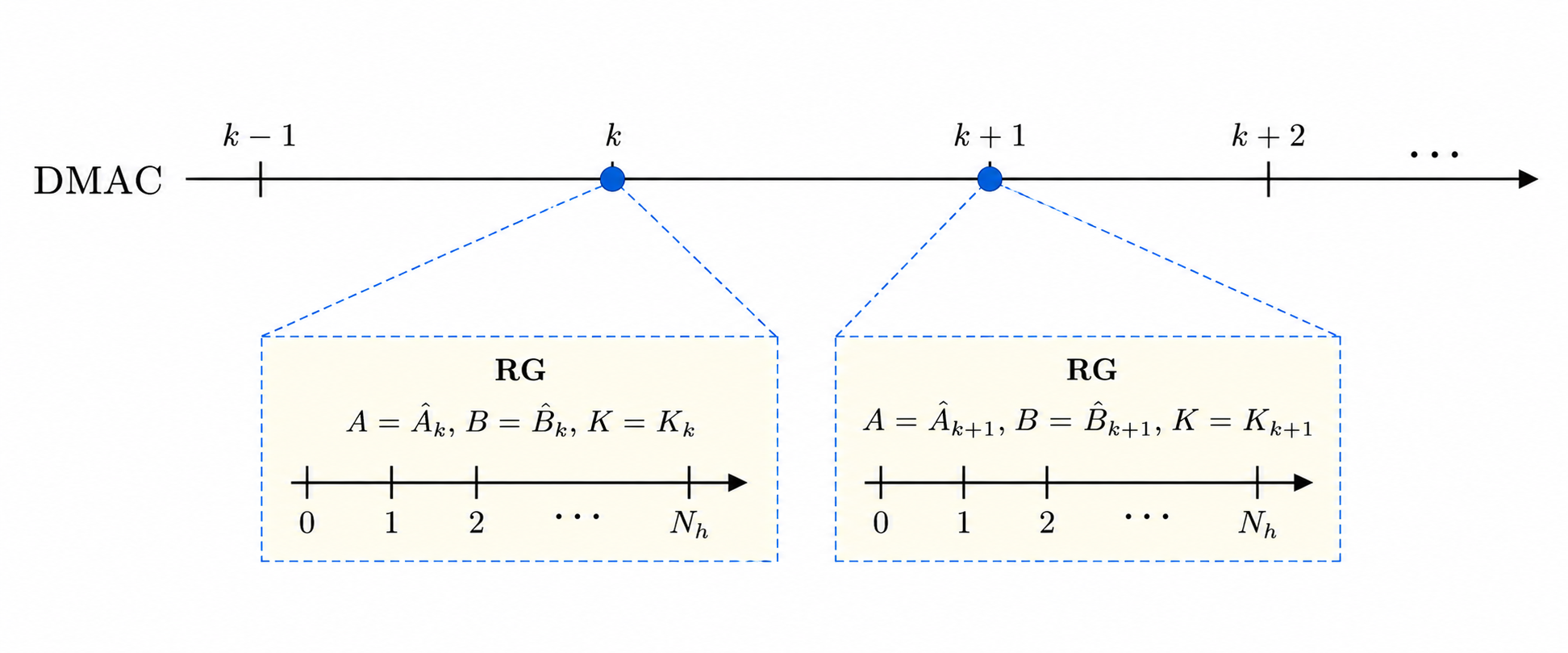}
    \caption{Frozen-model RG prediction at successive DMAC samples.}
    \label{fig:RG_frozen_prediction}
\end{figure}

As shown in Fig.~\ref{fig:RG_frozen_prediction}, the identified model and
feedback gains available at the current DMAC sample are held fixed during
the RG prediction.
At the next sample, the updated DMAC quantities define a new frozen
prediction model.

\subsection{Feasible-Reference Search}

Let the constrained output and input quantities be denoted by
$y_{\rmc,i|k}$ and $u_{\rmc,i|k}$, respectively, and suppose that they
must satisfy
\begin{align}
    y_{\rmc,\min}
    &\leq
    y_{\rmc,i|k}
    \leq
    y_{\rmc,\max},
    \label{eq:rg_output_constraints}
    \\
    u_{\rmc,\min}
    &\leq
    u_{\rmc,i|k}
    \leq
    u_{\rmc,\max}.
    \label{eq:rg_input_constraints}
\end{align}

For a candidate reference $g\in[g_{\min},g_{\max}]$, define the feasible
reference set
\begin{align}
    \mathcal{G}_k
    \isdef
    \Big\{
        g\in[g_{\min},g_{\max}]
        \;:\;&
        \eqref{eq:rg_output_constraints}
        \text{ and }
        \eqref{eq:rg_input_constraints}
        \text{ hold}
        \nonumber\\
        &\text{for }
        i=1,\ldots,N_h
    \Big\}.
    \label{eq:rg_feasible_set}
\end{align}

To balance tracking of the desired reference with smoothing of the
governed reference, define
\begin{align}
    J_{\rm ref}(g)
    =
    (g-r_k)^2
    +
    \lambda
    (g-g_{k-1})^2,
    \label{eq:rg_reference_cost}
\end{align}
where $\lambda\geq0$ weights the penalty on changes in the governed
reference.
In the examples considered in this paper, we set $\lambda=0.$

When $\mathcal{G}_k$ is nonempty, the governed reference is selected as
\begin{align}
    g_k
    =
    \arg\min_{g\in\mathcal G_k}
    J_{\rm ref}(g).
    \label{eq:rg_feasible_optimization}
\end{align}
The first term in \eqref{eq:rg_reference_cost} favors the desired
reference, while the second discourages large changes from the previously
applied governed reference.
Unlike the classical interpolation-based RG, this formulation does not
restrict $g_k$ to the segment between $g_{k-1}$ and $r_k$.

\begin{remark}
    The bounds $g_{\min}$ and $g_{\max}$ define the allowable reference-search
interval.
For the numerical examples considered in this paper, the upper search
bound is chosen as
    $g_{\max}
    =
    r_k,$
so that the governor does not search above the desired reference.
The lower bound $g_{\min}$ allows the governed reference to move below
both the desired and previously applied references when required by the
predicted constraints.
\end{remark}

\subsection{Reference Selection for an Empty Feasible Set}

The proposed method does not assume that the finite-horizon feasible set
$\mathcal G_k$ is always nonempty.
Changes in the identified prediction model, model error, disturbances, or
the current closed-loop state may result in $\mathcal G_k=\varnothing$.
In this case, no candidate reference satisfies all predicted constraints
over the specified horizon.
Rather than retaining $g_{k-1}$, which may itself produce predicted
constraint violations under the updated model, the RG selects the
candidate that minimizes the predicted constraint violation.

Let
$\mathcal Y_{\rmc}\isdef
[y_{\rmc,\min},y_{\rmc,\max}]$
and
$\mathcal U_{\rmc}\isdef
[u_{\rmc,\min},u_{\rmc,\max}]$
denote the admissible output and input sets, respectively.
The accumulated predicted constraint violation is defined as
\begin{align}
    J_{\rm viol}(g)
    =
    \sum_{i=1}^{N_h}
    \Big[
        &
        \operatorname{dist}^2
        \left(
            y_{\rmc,i|k}(g),
            \mathcal Y_{\rmc}
        \right)
        \nonumber\\
        &+
        \operatorname{dist}^2
        \left(
            u_{\rmc,i|k}(g),
            \mathcal U_{\rmc}
        \right)
    \Big],
    \label{eq:RG_violation_cost}
\end{align}
where $\operatorname{dist}(a,\mathcal A)$ denotes the Euclidean distance
from $a$ to the set $\mathcal A$.

When $\mathcal G_k=\varnothing$, the governed reference is selected as
\begin{align}
    g_k
    =
    \arg\min_{g\in[g_{\min},g_{\max}]}
    \left[
        J_{\rm viol}(g)
        +
        \rho J_{\rm ref}(g)
    \right],
    \label{eq:RG_infeasible}
\end{align}
where $\rho>0$ weights the reference-selection cost relative to the
predicted constraint-violation cost.
Thus, the dominant objective is to reduce predicted constraint violation,
while $J_{\rm ref}(g)$ retains a preference for tracking the desired
reference and avoiding unnecessary changes in the governed reference.

The current measured state is not included in $J_{\rm viol}$, only future
predicted violations are penalized.
Therefore, when satisfaction of all predicted constraints is impossible,
the RG selects the reference that minimizes the predicted violation over
the horizon.

The complete reference-selection rule is
\begin{equation}
\resizebox{0.882\columnwidth}{!}{$
g_k =
\begin{cases}
\displaystyle
\arg\min_{g\in\mathcal G_k}
J_{\rm ref}(g),
& \mathcal G_k\neq\varnothing,
\\[2mm]
\displaystyle
\arg\min_{g\in[g_{\min},g_{\max}]}
\left[
J_{\rm viol}(g)
+
\rho J_{\rm ref}(g)
\right],
& \mathcal G_k=\varnothing.
\end{cases}
$}
\label{eq:complete_rg_algorithm}
\end{equation}


\subsection{Numerical Implementation}
\label{sec:RG_numerical}

The reference-selection problems in
\eqref{eq:rg_feasible_optimization} and
\eqref{eq:RG_infeasible} are solved using a derivative-free direct-search
procedure based on a coarse grid with local refinement
\cite{conn2009introduction}.
At each sample, a finite set of candidate references is evaluated over
$[g_{\min},g_{\max}]$ using the frozen prediction model in
\eqref{eq:rg_augmented_prediction}.
If feasible candidates exist, the candidate minimizing
$J_{\rm ref}(g)$ is selected.
Otherwise, the candidate minimizing
$J_{\rm viol}(g)+\rho J_{\rm ref}(g)$ is selected.

After the best candidate is identified, the search interval is reduced
around that candidate and the direct search is repeated.
This refinement is performed for a prescribed number of levels to increase
the resolution of the selected governed reference.

\section{NUMERICAL RESULTS}
\label{sec:numerical_results}

The proposed DMAC--RG framework is evaluated on three examples of increasing complexity, including a coupled thermal system, a nonlinear Van der Pol oscillator, and a high-dimensional Burgers equation under sparse measurements.
%

\subsubsection{Linear System}

Consider the linear system
\begin{align}
    \dot{x}
    &=
    \begin{bmatrix}
        -0.70 & 0.30\\
         0.20 & -0.50
    \end{bmatrix}
    x
    +
    \begin{bmatrix}
        1.00\\
        0.20
    \end{bmatrix}u,
    \label{eq:lin_x}
    \\
    y
        &=
        \matl 
            0.6 & 0.4
        \matr
        x.
    \label{eq:lin_y}
\end{align}
The objective is to follow the command while constraining $y_{\rmc}=x_1.$ 
The state $x(t)$ is sampled at the period $T_\rms$, at which DMAC updates $u_k$ and the RG updates the governed reference $g_k$.
The numerical parameters of the system, DMAC, and the proposed RG are shown in Table \ref{tab:numerical_parameters}.

\begin{table}[ht]
\centering
\caption{Simulation, DMAC, and RG parameters used in the numerical examples.}
\label{tab:numerical_parameters}
\scriptsize
\renewcommand{\arraystretch}{1.08}
\setlength{\tabcolsep}{2.2pt}

\begin{tabular}{p{0.27\columnwidth}
                p{0.21\columnwidth}
                p{0.21\columnwidth}
                p{0.21\columnwidth}}
\hline
\textbf{Parameter} &
\centering\textbf{Linear System} &
\centering\textbf{Van der Pol} &
\centering\arraybackslash\textbf{Burgers} \\
\hline

\multicolumn{4}{l}{\textbf{Plant Parameters}} \\
\hline

Plant dimension
& \centering $2$
& \centering $2$
& \centering\arraybackslash $100$ \\

Sampling period $T_s$
& \centering $0.1$
& \centering $0.1$
& \centering\arraybackslash $0.01$ \\

Constraint $y_{\rmc,\max}$
& \centering $1$
& \centering $1$
& \centering\arraybackslash $0.7$ \\

\hline

\multicolumn{4}{l}{\textbf{DMAC Parameters}} \\
\hline

$Q_{\mathrm{lqr}}$
& \centering $I_3$
& \centering $I_3$
& \centering\arraybackslash $10I_8$ \\

$R_{\mathrm{lqr}}$
& \centering $1$
& \centering $1$
& \centering\arraybackslash $0.1$ \\

DMAC $\lambda$
& \centering $0.995$
& \centering $0.995$
& \centering\arraybackslash $0.9995$ \\

DMAC $R_\Theta$
& \centering $10I_3$
& \centering $10I_3$
& \centering\arraybackslash $100I_8$ \\

Excitation noise
& \centering $10^{-3}\mathcal N(0,1)$
& \centering $10^{-3}\mathcal N(0,1)$
& \centering\arraybackslash $10^{-3}\mathcal N(0,1)$ \\

\hline

\multicolumn{4}{l}{\textbf{RG Parameters}} \\
\hline

Horizon $N_h$
& \centering $200$
& \centering $200$
& \centering\arraybackslash $200$ \\

Search interval $g$
& \centering $[-2,r]$
& \centering $[-2,r]$
& \centering\arraybackslash $[-2,r]$ \\

Initial grid size $N_g$
& \centering $21$
& \centering $21$
& \centering\arraybackslash $21$ \\

Refinement levels
& \centering $2$
& \centering $2$
& \centering\arraybackslash $2$ \\

Weight $\rho$
& \centering $10^{-3}$
& \centering $10^{-3}$
& \centering\arraybackslash $10^{-3}$ \\

\hline
\end{tabular}
\end{table}

%

\begin{figure}[ht]
\centering
\includegraphics[width=\linewidth]{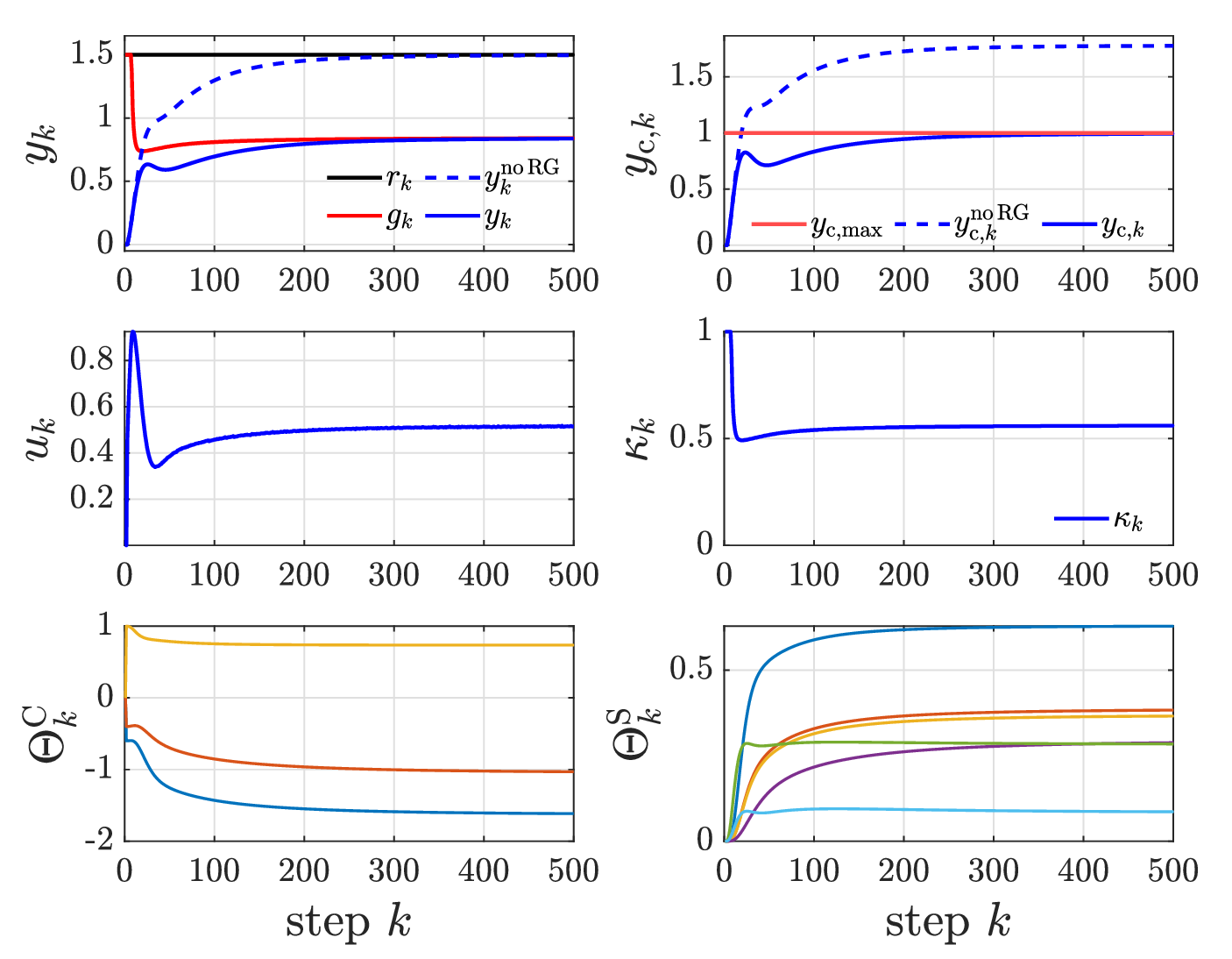}
\caption{DMAC--RG response for the linear system, with and without the reference governor.
The first row shows the regulated output $y_k$ and constrained output $y_{\rmc,k}$,
the second row shows the control signal $u_k$ and interpolation factor $\kappa_k$, and 
the third row shows the DMAC controller gain $\Theta_k^\rmC$ and estimated plant parameters $\Theta_k^\rmS$.}
\label{fig:DMAC_RG_Linear_V2}
\end{figure}

Figure \ref{fig:DMAC_RG_Linear_V2} shows the closed-loop response of the system \eqref{eq:lin_x}-\eqref{eq:lin_y} with DMAC in the loop and the proposed reference governor modifying the commanded reference.
The figure shows the outputs $y_k$ and $y_{\rmc,k}$ both with and without the reference governor.
With the reference governor, $y_{\rmc,k}$ approaches $y_{\rmc,\max}$ and $y_k$ converges to the largest value consistent with $y_{\rmc,k}\leq y_{\rmc,\max}$.
Note that the proposed RG does not use the true plant model $(A,B)$ to compute $g_k$.
Instead, it uses the estimates provided by DMAC to synthesize the governed reference.

The second and third rows in Figure \ref{fig:DMAC_RG_Linear_V2} show, respectively, the control signal $u_k$ computed by DMAC together with the equivalent interpolation factor
\begin{align}
    \kappa_k
    \isdef
    \frac{g_k-g_0}{r_k-g_0},
    \label{eq:kappa_k_def}
\end{align}
and the controller gain $\Theta_k^\rmC$ and estimated plant parameters $\Theta_k^\rmS$ computed by the DMAC algorithm.
Unlike the classical RG in Section~\ref{sec:classicalRG}, $\kappa_k$ here is not an optimization variable, it is reported only as a diagnostic quantity for comparison.


\subsubsection{Van der Pol Oscillator}

Consider the nonlinear Van der Pol oscillator
\begin{align}
    \dot{x}_1
    &=
    x_2,
    \label{eq:vdp_x1}
    \\
    \dot{x}_2
    &=
    \mu(1-x_1^2)x_2
    -
    x_1
    +
    u,
    \label{eq:vdp_x2}
    \\
    y
    &=
    \matl
        0.6 & 0.4
    \matr
    x,
    \label{eq:vdp_y}
\end{align}
where $\mu=1$.
The objective is to follow the command while constraining $y_{\rm c}=x_1$.
The state $x(t)$ is sampled at the period $T_{\rm s}$, at which DMAC updates $u_k$ and the RG updates the governed reference $g_k$.
The numerical parameters of the system, DMAC, and the proposed RG are shown in Table~\ref{tab:numerical_parameters}.

\begin{figure}[ht]
\centering
\includegraphics[width=\linewidth]{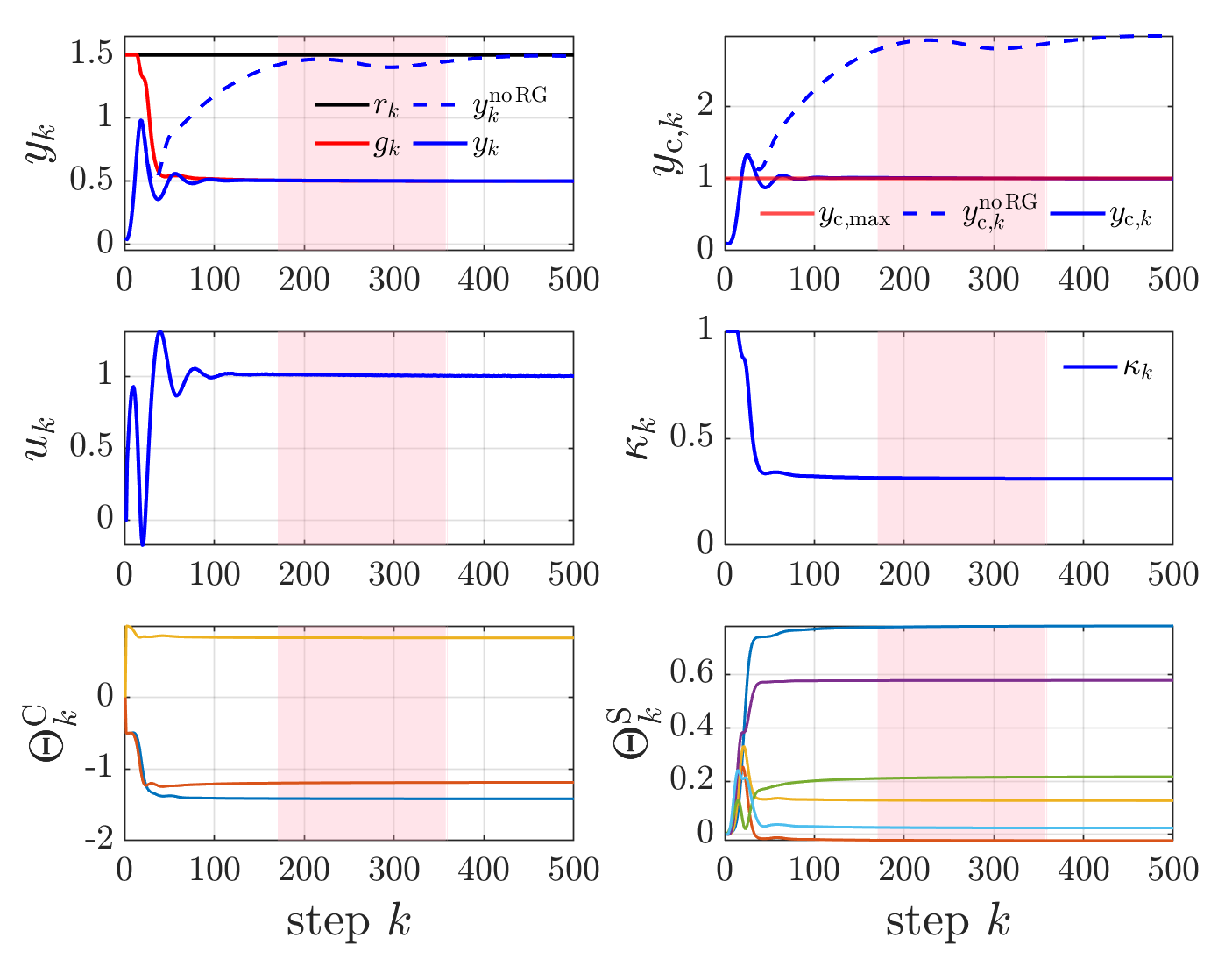}
\caption{DMAC--RG response for the nonlinear Van der Pol oscillator, with and without the reference governor. The first row shows the regulated output $y_k$ and constrained output $y_{{\rm c},k}$, the second row shows the control signal $u_k$ and interpolation factor $\kappa_k$, and the third row shows the DMAC controller gain $\Theta_k^{\rm C}$ and estimated plant parameters $\Theta_k^{\rm S}$.}
\label{fig:DMAC_RG_VDP_V2}
\end{figure}

Figure~\ref{fig:DMAC_RG_VDP_V2} shows the closed-loop response of the system \eqref{eq:vdp_x1}--\eqref{eq:vdp_y} with DMAC in the loop and the proposed reference governor modifying the commanded reference.
The pink shaded region marks the period during which $\SG_k = \varnothing$ and thus $g_k$ is given by \eqref{eq:RG_infeasible}.
The figure shows the outputs $y_k$ and $y_{{\rm c},k}$ both with and without the reference governor.
With the reference governor, $y_{{\rm c},k}$ approaches $y_{{\rm c},\max}$ while satisfying the constraint.
Note that the proposed RG does not use the nonlinear Van der Pol model to compute $g_k$. 
Instead, it uses the model identified online by DMAC in the reference governor.
The second and third rows in Figure~\ref{fig:DMAC_RG_VDP_V2} show, respectively, the control signal $u_k$ computed by DMAC together with the equivalent interpolation factor defined in \eqref{eq:kappa_k_def} and the controller gain $\Theta_k^{\rm C}$ and estimated plant parameters $\Theta_k^{\rm S}$ computed by the DMAC algorithm.
%

\subsubsection{Burgers Equation with Sparse Measurements}

Consider the nonlinear viscous Burgers equation
\begin{align}
    \frac{\partial w}{\partial t}
    +
    w\frac{\partial w}{\partial s}
    &=
    \nu\frac{\partial^2 w}{\partial s^2}
    +
    b(s)u(t),
    \label{eq:burgers}
\end{align}
where $\nu=0.1$.
To simulate \eqref{eq:burgers}, we spatially discretize the $s-$domain using $100$ nodes to obtain an ordinary differential equation.
The discretized state $w(t)$ of the resulting ODE is sampled at the period $T_\rms$, at which DMAC updates $u_k$ and the RG updates the governed reference $g_k$.
%
In DMAC, the partial state is
\begin{align}
    \xi_k
    =
    \big[
        &x_{1,k},\;
        x_{16,k},\;
        x_{31,k},\;
        x_{46,k}, \nonumber\\[-1mm]
        &x_{61,k},\;
        x_{76,k},\;
        x_{91,k}
    \big]^{\top}.
    \label{eq:burgers_sparse_state}
\end{align}
The regulated output is
\begin{align}
    y_k = x_{61,k},
    \label{eq:burgers_y}
\end{align}
and the objective is to follow the command while constraining
$y_{{\rm c},k}=x_{1,k}$.
The control input acts at grid node $55$. 
The numerical parameters of the system, DMAC, and the proposed RG are shown in Table~\ref{tab:numerical_parameters}.

\begin{figure}[ht]
\centering
\includegraphics[width=\linewidth]
{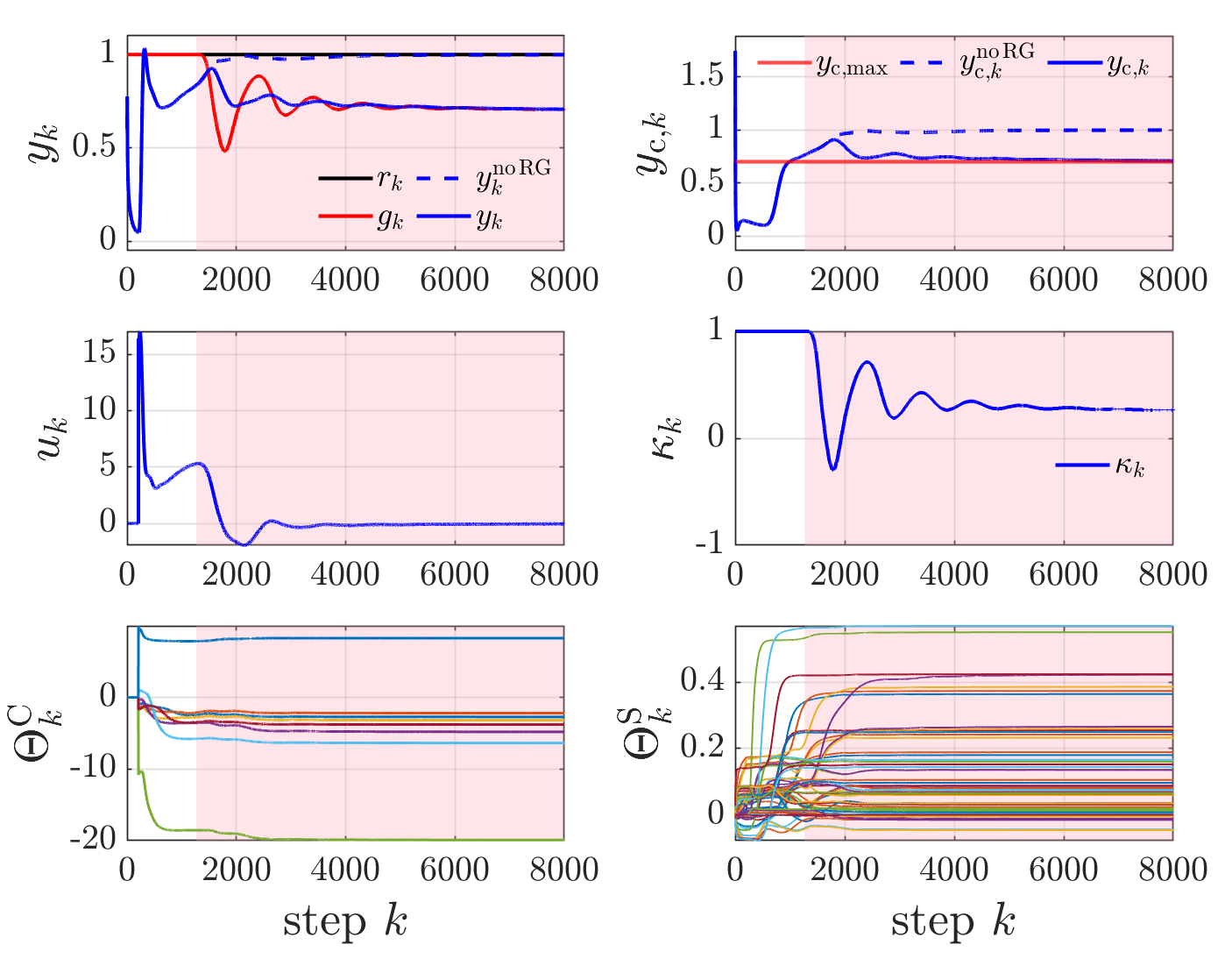}
\caption{DMAC--RG response for the Burgers equation using seven measured states, with and without the reference governor. The first row shows the regulated output $y_k$ and constrained output $y_{{\rm c},k}$, the second row shows the control signal $u_k$ and interpolation factor $\kappa_k$, and the third row shows the DMAC controller gain $\Theta_k^{\rm C}$ and estimated plant parameters $\Theta_k^{\rm S}$.}
\label{fig:DMAC_RG_Burgers_PartialStates_V2}
\end{figure}

\begin{figure}[ht]
\centering
\includegraphics[width=\linewidth]
{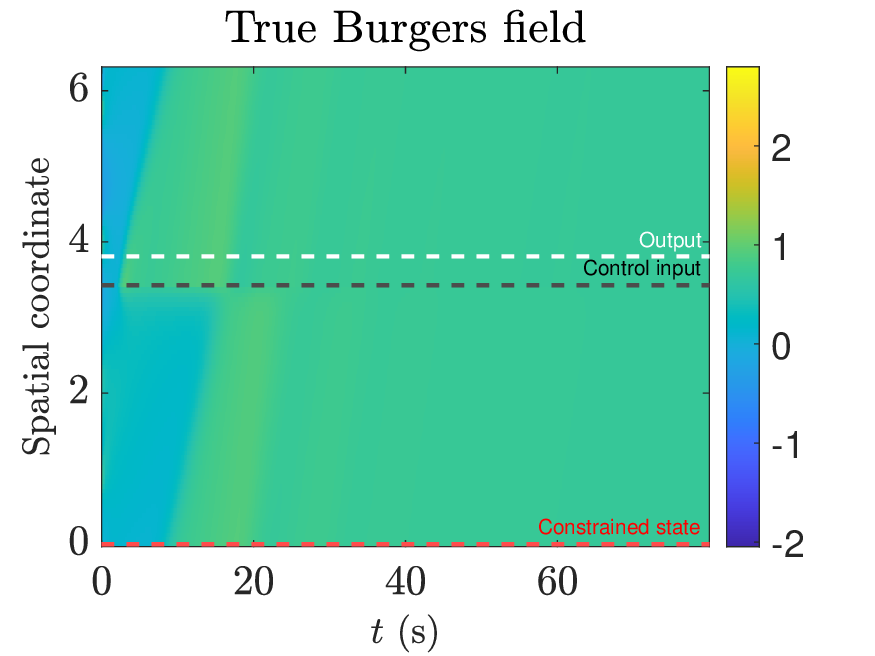}
\caption{Spatiotemporal evolution of the $100$-state Burgers field. The dashed lines indicate the spatial locations of the regulated output, constrained output, and control input.}
\label{fig:DMAC_RG_Burgers_PartialStates_Contour_V2}
\end{figure}

Figure~\ref{fig:DMAC_RG_Burgers_PartialStates_V2} shows the closed-loop response of the system \eqref{eq:burgers} with DMAC in the loop and the proposed reference governor modifying the commanded reference.
%
The pink shaded region marks the period during which $\SG_k = \varnothing$ and thus $g_k$ is given by \eqref{eq:RG_infeasible}.
With the reference governor, $y_{{\rm c},k}$ approaches $y_{{\rm c},\max}$ while satisfying the constraint.
Note that the proposed RG does not use the complete $100$-state Burgers model to compute $g_k$. Instead, it uses the seventh-order model identified online by DMAC from the sparse measurements.
The second and third rows in
Figure~\ref{fig:DMAC_RG_Burgers_PartialStates_V2} show, respectively, the
control signal $u_k$ computed by DMAC together with the equivalent
interpolation factor defined in \eqref{eq:kappa_k_def}
and the controller gain $\Theta_k^{\rm C}$ and estimated plant parameters $\Theta_k^{\rm S}$ computed by the DMAC algorithm.

In this example, DMAC and the RG operate using only the seven measured states, rather than the complete $100$-state field.

Figure~\ref{fig:DMAC_RG_Burgers_PartialStates_Contour_V2} shows the corresponding evolution of the complete spatial field. 
Although the complete field is shown for visualization, it is not used by either DMAC or the RG. 
Only the seven measured states in \eqref{eq:burgers_sparse_state} are used for online identification, control, and RG prediction.

\section{CONCLUSION}
\label{sec:conclusion}
This paper developed a reference-governor framework for constraint enforcement in Dynamic Mode Adaptive Control.
The RG predicts using the model identified online, enabling constraint enforcement without a prior plant model and with partial-state measurements.
Classical interpolation-based RG assumes a fixed reference-governable set; the proposed method instead searches directly over the governed reference to accommodate an online-changing model.
When a feasible reference exists over the prediction horizon, the reference closest to the desired command is selected; when none exists, the reference minimizing predicted constraint violation is selected, with a secondary preference for reference tracking.
Numerical examples on linear and nonlinear systems, including a high-dimensional system under partial-state measurement, demonstrated constraint enforcement using online-identified dynamics in place of the true plant model.

Future work will develop conditions under which feasibility can be guaranteed for the reference governor under general indirect adaptive control, as the identified model evolves online.

\printbibliography
\end{document}